\documentclass[12pt,thmsa]{article}
\usepackage{amsmath, latexsym, amsfonts, amssymb, amsthm, amscd}
\usepackage{graphics}
\newtheorem{theorem}{Theorem}

\newtheorem{proposition}{Proposition}
\newtheorem{lemma}{Lemma}
\newtheorem{definition}{Definition}
\newtheorem{rem}{Remark}

\newcommand{\p}{\Bbb{P}}

\newcommand{\e}{\Bbb{E}}

\newcommand{\ind}{\mbox{\rm 1\hspace{-0.04in}I}}

\newcommand{\R}{\mbox{\rm I\hspace{-0.02in}R}}

\newcommand{\el}{\stackrel{\mathcal L}{=}}
\newcommand{\nel}{\stackrel{\mathcal L}{\neq}}

\newcommand{\eqdef}{\stackrel{\mbox{\tiny$($def$)$}}{=}}

\def\QED{\hfill\vrule height 1.5ex width 1.4ex depth -.1ex \vskip20pt}

\begin{document}

\title{Reflection principle and Ocone martingales.}

\maketitle

\begin{center}
{\large L. Chaumont}\footnote{$^{,2}$ LAREMA, D\'epartement de
Math\'ematiques, Universit\'e d'Angers, 2, Bd Lavoisier - 49045,
\\\hspace*{.4in}{\sc Angers Cedex 01.}

\hspace*{.05in}$^1$E-mail: loic.chaumont@univ-angers.fr$\;\;\;$
$^2$E-mail: lioudmila.vostrikova@univ-angers.fr}{\large
 and  L. Vostrikova$^2$}
\end{center}
\vspace{0.2in}

\begin{abstract} Let $M =(M_t)_{t\geq 0}$ be any continuous real-valued stochastic
process. We prove that if there exists a sequence $(a_n)_{n\geq 1}$ of real
numbers which converges to 0 and such that $M$ satisfies the
reflection property at all levels $a_n$ and $2a_n$ with $n\geq 1$, then $M$ is an
Ocone local martingale with respect to its natural filtration.
We state the subsequent open question: is this result still true
when the property only holds at levels $a_n$~? Then we prove that the later question is
equivalent to the fact that for Brownian motion, the $\sigma$-field of the invariant events
by all  reflections at levels $a_n$, $n\ge1$ is trivial.
We establish similar results for skip free $\mathbb{Z}$-valued processes and use them
for the proof in continuous time,  via a discretisation in space.\\


\end{abstract}

\section{Introduction and main results}

Local martingales whose law is invariant under any integral transformations preserving their
quadratic variation were first introduced and characterized by Ocone \cite{oc}.
Namely a continuous real-valued local martingale
$M=(M_t)_{t\geq 0}$ with natural filtration
$\mathbb F = ({\cal F}_t)_{t\geq 0}$ is called {\it Ocone}
if
\begin{equation}\label{a}
\left(\int_0^t H_s dM_s\right)_{t\geq 0}
\stackrel{\mathcal L}{=} M\,,
\end{equation}
for all processes $H$ belonging to the set
$$\mathcal H = \{ H=(H_t)_{t\geq 0}\,\,|\,H \mbox{\,is\, $\mathbb F$-predictable},\, |H_t|
=1,\, \mbox{\,for \,all\,\,} t\geq 0\}.$$

In the primary paper \cite{oc}, the author
proved that a local martingale is Ocone whenever it satisfies
(\ref{a}) for all processes $H$ belonging to the smaller class
of deterministic processes:
\begin{equation}
\label{b} \mathcal H_1=\{\,\left(\ind _{[0,u]}(t) - \ind
_{]u,+\infty [}(t)\right)_{t\geq 0}\,,\,\mbox{with\,\,}u\geq 0\}\,.
\end{equation}
A natural question for which we sketch out an answer in this paper
is to describe minimal sub-classes of $\mathcal H$ characterizing
Ocone local martingales through relation (\ref{a}). For instance, it is
readily seen that the subset $\{\,\left(\ind _{[0,u]}(t) - \ind
_{]u,+\infty [}(t)\right)_{t\geq 0}\,,\,\mbox{with\,\,}u\in E\}$ of
$\mathcal H_1$ characterizes Ocone martingales if and only if $E$ is
dense in $[0,\infty)$. Let us denote by
$\langle M\rangle $ the quadratic variation of $M$. In \cite{oc} it was shown
that for continuous local martingales, (\ref{a}) is equivalent to the fact
that conditionally to the $\sigma $-algebra $\sigma\{ \langle M
\rangle _s, s \geq 0 \}$, $M$ is a gaussian process with independent
increments. Hence {\it a continuous Ocone local martingale is a Brownian motion time changed
by any independent nondecreasing continuous process.} This is actually the definition
we will refer to all along this paper.

When the continuous local martingale $M$ is divergent, i.e. $\p$-a.s.
\[\lim_{t\rightarrow\infty}\langle M\rangle_t = +\infty\,,\]
we denote by $\tau$ the right continuous inverse of $\langle
M\rangle$, i.e. for $t\ge0$,
$$\tau_t=\inf\{s\ge0:\langle M\rangle_s>t\}\,,$$
and we recall that the Dambis-Dubins-Schwarz Brownian motion
associated to $M$ is the $({\cal F}_{\tau_t})$-Brownian motion
defined by
$$B^M\eqdef (M_{\tau_t})_{t\ge0}.$$
Then Dubins, Emery and Yor \cite{dey} refined Ocone's characterization
by proving that (\ref{a}) is equivalent to each of the following
three properties:
\begin{itemize}
\item[(i)] The processes $\langle M\rangle$ and
$B^M$ are independent.
\item[(ii)] For every $\mathbb F$-predictable process $H$, measurable
for the product $\sigma$-field ${\cal B}(\mathbb{R}_+)\otimes\sigma(\langle
M\rangle)$ and such that $\int_0^\infty H^2_s\,d\langle
M\rangle_s<\infty$, $\p$-a.s.,
\[\e\left(\exp\left(i\int_0^\infty H_s\,dM_s\right)\,|\,\langle
M\rangle\right)=\exp\left(-\frac12\int_0^\infty H^2_s\,d\langle
M\rangle_s\right)\,.\]
\item[(iii)] For every deterministic function $h$ of the form
$\sum_{j=1}^n\lambda_j\ind_{[0,a_j]}$,
\[\e\left[\exp\left(i\int_0^\infty h(s)\,dM_s\right)\right]=
\e\left[\exp\left(-\frac12\int_0^\infty h^2(s)\,d\langle
M\rangle_s\right)\right]\,.\]
\end{itemize}
It can be easily shown that the equivalence between (\ref{a}) and (i), (ii), (iii) also holds
in the case when $M$ is not necessarily
divergent. This fact will be used in the proof of Theorem \ref{main1}.
We also refer to \cite{vy} for further results
related to Girsanov theorem and different classes of martingales.

In  \cite{dey}, the authors conjectured that the class $\mathcal H_1$ can be reduced to
a single process, namely that (\ref{a}) is equivalent to:
\begin{equation}\label{sgn}
\left(\,\int_0^{t}\mbox{sign}(M_s)\,dM_s\,\right)_{\,t\ge0}\el M\,.
\end{equation}
In fact, (\ref{sgn}) holds if and only if $B^M$ and
$\langle M \rangle $ are conditionally independent given
 the $\sigma$-field of
invariant sets by the L\'evy transform of $B^M$, i.e.
$B^M\mapsto \left(\int_0^{\cdot}\mbox{sign}(B^M_s)\,dB^M_s\right)$, see \cite{dey}.
Hence if the L\'evy transform of Brownian motion is ergodic, then
 $B^M$ and $\langle M \rangle$ are independent and (\ref{sgn})
implies that $M$ is an Ocone local martingale. The converse is also proved in \cite{dey}, that is
if (\ref{sgn}) implies that $M$ is an Ocone local martingale, then the L\'evy transform of
Brownian motion is ergodic.

Different other approaches have been proposed to prove
ergodicity of the L\'evy transform but this problem is still open.
Among the most accomplished works in this direction, we may cite
papers by Malric \cite{ma1}, \cite{ma2} who studied the density of
zeros of iterated L\'evy transform. Let us also mention
that in discrete time case  this problem has been
treated in \cite{ds} where the authors proved that an equivalent of the L\'evy transform for
symmetric Bernoulli random walk is ergodic.\\

In this paper we exhibit a new  sub-class of $\mathcal H_1$
characterizing continuous Ocone local martingales which is related to
first passage times and the reflection property of stochastic
processes. If $M$ is the standard Brownian motion and $T_a(M)$ the
first passage time at level $a$, i.e.
\begin{equation}\label{ta}
T_a(M)=\inf\{t\geq 0:M_t=a\},
\end{equation}
where here and in all the remainder of this article, we make the convention that
$\inf \{\emptyset\}= +\infty$, then for all $a\in\mathbb{R}$:
$$(M_t)_{t\geq 0}\el ( M_t\ind _{\{ t\leq T_a(M)\}} +
(2a-M_t)\ind _{\{t > T_a(M)\}})_{t\geq 0}.$$
It is readily checked that this identity in law actually holds for any
continuous Ocone local martingale. This property  is
known as the {\it reflection principle at level $a$} and was first observed
for symmetric Bernoulli random  walks by Andr\'e \cite{an}. We will use
this terminology for any continuous stochastic process $M$ and when
no confusion is possible, we will denote by $T_a=T_a(M)$ the first passage
time at level $a$ by $M$ defined as above.

Let $(\Omega ,\mathcal F , \mathbb F, \p)$ be the canonical space of
continuous functions endowed with its natural right-continuous filtration $\mathbb F =
(\mathcal F _t)_{t\geq 0}$ completed by negligible sets of $\mathcal
F= \bigvee_{t\geq 0} \mathcal F_t$.
The family of transformations $\Theta^a$, $a\ge0$, is
defined for all continuous functions $\omega\in\Omega$ by
\begin{equation}\label{transf}
\Theta^{a}(\omega)=(\omega_t\ind _{\{ t\leq T_a\}} +
(2a-\omega_t)\ind _{\{t > T_a\}})_{t\geq 0}\,.
\end{equation}
 Note that
$\Theta^a(\omega)=\omega$ on the set $\{\omega : T_a(\omega)=\infty\}$. When
$M$ is a local martingale, $\Theta^{a}(M)$ can by expressed in terms
of a stochastic integral, i.e.
\[\Theta^{a}(M)=\left(\int_0^t\,\left(\ind _{[0,T_a]}(s) - \ind
_{]T_a,+\infty [}(s)\,\right)\,dM_s\right)_{t\ge0}.\]
The set
$\mathcal H_2 = \{\left(\ind _{[0,T_a]}(t) -
\ind_{]T_a,+\infty[}(t)\right)_{t\ge0}\,|\,a\ge0\}$ is a subclass of
$\mathcal H$ which provides a family of transformations preserving
the quadratic variation of $M$ and we will prove that it characterizes Ocone local
martingales. But the fact that the transformations $\omega\mapsto \Theta^{a}(\omega)$ are defined
for all continuous functions  $\omega\in\Omega$ allows us to
characterize Ocone local martingales in the whole set of continuous stochastic
processes as shows our main result.

\begin{theorem}\label{main1}
Let $M =(M_t)_{t\geq 0}$ be a continuous stochastic process defined on the canonical probability
space, such that $M_0=0$.
If there exists a sequence $(a_n)_{n\geq 1}$ of positive real numbers
such that $\lim_{n\rightarrow \infty} a_n=0$ and for all
$n\ge0$:
\begin{equation}\label{hyp}
\Theta^{a_n}(M)\el \Theta^{2a_n}(M)\el M\,,
\end{equation}
then $M$ is an Ocone local  martingale with
respect to its natural filtration. Moreover, if $T_{a_1}<\infty$
a.s., then $M$ is a divergent local martingale.
\end{theorem}

\begin{rem} It is natural to wonder about the necessity of the
hypothesis $\Theta^{2a_{n}}(M)\el M$ in Theorem $\ref{main1}$. The discrete time counterpart
of this problem which is presented in section $\ref{skipfree}$, shows that it is necessary
for a skip free process $M$ to satisfy $\Theta^a(M)\el M$, for $a=0,1$
and $2$ in order to be a skip free   Ocone local martingale, i.e. the reflection property
at $a=0$ and $1$ is not sufficient, see the counterexamples in section $\ref{counterexample}$.
This argument seems to confirm that the assumption  $\Theta^{2a_{n}}(M)\el M$ is necessary
 in continuous time.
\end{rem}

In an attempt to identify the sequences $(a_n)_{n\geq 1}$ which characterize Ocone local martingales,
 we obtained the following theorem. Let
$a=(a_n)_{n\geq 1}$ be a sequence of real numbers with
$\lim_{n\rightarrow\infty} a_n=0$ and
let ${\mathcal I}^a$ the sub-$\sigma$-field  of
 the invariant sets by all the transformations $\Theta^{a_n}$, i.e.
\[{\cal I}^{a}=\{F\in{\cal F}:\ind_{F}\circ\Theta^{a_n}\stackrel{\mbox{\small a.s.}}{=}
 \ind_F,\;\;\mbox{for all $n\ge0$}\}.\]

\begin{theorem}\label{main2}
The following assertions are equivalent:
\begin{itemize}
\item[$(i)$]  Any continuous  local martingale $M$ satisfying
$\Theta^{a_n}(M)\el M$ for all $n\ge0$ is an Ocone local martingale.
\item[$(ii)$] The sub $\sigma$-field
 ${\cal I}^{a}$ is trivial for the Wiener measure on the canonical space $(\Omega ,\mathcal F)$.
\end{itemize}
\end{theorem}

\begin{rem} It follows from Theorems $\ref{main1}$ and $\ref{main2}$ that if the sequence
$(a_n)$ contains a subsequence $(2a_{n'})$ $( $this holds, for instance, when $(a_n)$
is dyadic sequence $)$, then the sub $\sigma$-field
${\cal I}^{a}$ is trivial for the Wiener measure on $(\Omega,\mathcal F)$. So, our open question
is equivalent to: is the sub $\sigma$-field ${\cal I}^{a}$  trivial for any sequence $(a_n)$
decreasing to zero~?
\end{rem}

In the next section,  we prove analogous results for skip
free processes. We use them as preliminary results to prove Theorem
\ref{main1} in section \ref{proof}.  In section \ref{counterexample}, we give
 counterexamples in the discrete time setting, related to Theorem \ref{prop2}.
Finally, in section \ref{conj},  we prove Theorem \ref{main2}.

\section{Reflecting property and skip free processes}\label{skipfree}
\subsection{Discrete time skip free processes}
A discrete time skip free process $M$ is any measurable
stochastic process with $M_0=0$ and for all $n\ge1$,
$\Delta M_n= M_{n}-M_{n-1}\in\{-1,0,1\}$. This section is devoted to an analogue of
Theorem \ref{main1} for skip free processes.\\

To each skip free process $M$, we associate the increasing process
$$[M]_n = \sum_{k=0}^{n-1}(M_{k+1}-M_k)^2\,,\;\;n\ge1\,,\;\;\;[M]_0=0\,,$$
which is called the {\it quadratic variation} of $M$. In this section,
since no confusion is possible, we will use the same notations for discrete
processes as in continuous time case. For every integer $a\ge0$, we denote
by $T_a$ the first passage time by $M$ to the level $a$,
$$T_a=\inf\{k\geq 0:M_k=a\}\,.$$
We also introduce the inverse process $\tau$ which is defined by $\tau_0=0$ and for $n\geq 1$,
$$\tau_n=\inf\{k>\tau_{n-1}:[M]_k=n\}$$
with $\inf\{\emptyset\} = \tau_{n-1}$.
 Then we may
define
\begin{equation}\label{sup}
S^M=(M_{\tau_n})_{n\ge0}
\end{equation}
Denote also
$$T=\inf\{k\geq 0 \,:\, [S^ M]_k=[S^ M]_{\infty}\}\,,$$
then note that the paths of $S^M$ are such that $\Delta S^M_k\in\{-1,+1\}$, for all $k\leq T$
and $\Delta S_k^M=0$, for all  $k>T$, where  $\Delta S_k^M= S_k - S_{k-1}^M$.\\

We recall that skip free martingales are just skip free processes being martingales
with respect to some filtration.
It is well known that for any divergent free skip martingale $M$, that is satisfying
$\lim_{n\rightarrow+\infty}[M]_n=+\infty$, a.s., the process $S^M$ is a symmetric
Bernoulli random walk on $\mathbb Z$. This property is the equivalent of the
Dambis-Dubins-Schwartz theorem for continuous martingales. In discrete
time, the proof is quite straightforward and we recall it now.

A first step is the equivalent of L\'evy's characterization for skip free martingales~:
any skip free martingale $S$ such that $S_{n+1}-S_n\neq0$, for all
$n\ge0$ (or equivalently, whose quadratic variation satisfies
$[S]_n=n$) is a symmetric Bernoulli random walk. Indeed for $n\ge1$,
$S_{1},S_{2}-S_{1},\dots,S_{n}-S_{n-1}$ are i.i.d. symmetric Bernoulli
r.v.'s
if and only if for any subsequence $1\le n_1\le\dots\le n_k\le n$:
\begin{eqnarray*}
&&\e[(S_{n_1}-S_{n_1-1})(S_{n_2}-S_{n_2-1})\dots(S_{n_k}-S_{n_k-1})]=\\
&&\quad\qquad\e[S_{n_1}-S_{n_1-1}]\e[S_{n_2}-S_{n_2-1}]\dots\e[S_{n_k}-S_{n_k-1}]=0
\end{eqnarray*}
and this identity  can  be easily checked from the martingale property.
Finally call ${\mathbb F}=({\cal F}_n)_{n\ge0}$ the natural filtration generated by $M$.
Since $[M]_n$ is an ${\mathbb F}$-adapted process, from the optional stopping theorem,
$S^M$ is a martingale with respect to the filtration $({\cal
F}_{\tau_n})_{n\ge0}$ and since its increments cannot be $0$, we
conclude from L\'evy's characterization.

We  recall also  the following important property:
any skip free process which is  a symmetric Bernoulli random walk time
changed by an independent nondecreasing skip free process,
is a local martingale with respect to its natural filtration.\\

This leads to the definition:
\begin{definition}\label{def}
A discrete Ocone local martingale is a symmetric Bernoulli random walk time changed by any independent
increasing skip free process.
\end{definition}

We emphasize that in this particular case, Definition \ref{def}
coincides with the general definition of Ocone \cite{oc}.
It should also be noticed that the symmetric Bernoulli random walk of Definition 1
is not necessarily the same as in  (\ref{sup}). It coincides with $S^M$ if $M$ is a divergent process.
If $M$ is not divergent, then it can obtained obtained from the initial one by pasting of
an independent symmetric Bernoulli random walk (see Lemma 3), otherwise the independence can fail.
\vspace{0.2cm}

\noindent A counterpart of transformations $\Theta^a$ defined in (\ref{transf}) for skip free processes
is given for all integers $a\ge0$ by
\begin{equation}\label{trdisc}
\Theta^a(M)_n=\sum_{k=1}^n(\ind_{\{k \leq T_a\}}- \ind_{\{k >
T_a\}})\Delta M_{k}\,,
\end{equation}
where $\Delta M_k=M_k-M_{k-1}$. Again in the following discrete time counterpart
of Theorem~\ref{main1}, we characterize discrete Ocone local martingales in the whole set of skip free
processes.

\begin{theorem}\label{prop2}
Let $M$ be any discrete  skip free process. Assume that  for all $a\in\{0,1,2\}$,
\begin{equation}\label{identity}
\Theta^a(M)\el M,\end{equation}
then $M$ is a discrete Ocone local martingale with respect to its
natural filtration. If in addition $T_1 < \infty $ a.s. then $M$ is a
divergent local martingale.
\end{theorem}

\noindent The proof of Theorem \ref{prop2} is based on  the following
crucial combinatorial lemma concerning the set of sequences of
partial sums of elements in $\{-1,+1\}$ with length $m\ge1$:
\[\Lambda^m=\{(s_0,s_1,\dots,s_m):s_0=0\;\;\mbox{and}\;\; \Delta s_{k}
\in\{-1,+1\} \,\mbox{for}\,\;\; 1\leq k\leq m\,\},\]
where $\Delta s_k=s_k - s_{k-1}$.

For each sequence $s\in\Lambda^m$, and each integer $a$, we define
$T_a(s)=\inf\{k\geq 0:s_k=a\}$, with $\inf\emptyset=+\infty$. The transformation $\Theta^a(s)$ is defined
for each $s\in\Lambda^m$ by
\[\Theta^a(s)_n=\sum_{k=1}^n(\ind_{\{k\leq T_a(s)\}}- \ind_{\{k >
T_a(s)\}})\Delta s_{k}\,,\;\;\;n\le m.\]

\begin{lemma}\label{lem1}Let $m\geq1$ be fixed.
For any two elements $s$ and $s'$ of the set $\Lambda^m$ such that
$s\neq s'$, there are integers $a_1,a_2,\dots,a_k\in\{0,1,2\}$
depending on $s$ and $s'$ such that
\begin{equation}\label{ssp}
s'=\Theta^{a_k}\Theta^{a_{k-1}}\dots\Theta^{a_1}(s)\,.
\end{equation}
Moreover,  the integers $a_1,\dots,a_k$ can be chosen so that
$s\in\Lambda_{a_1}^m$ and $\Theta^{a_{i-1}}\Theta^{a_{i-2}}\dots\Theta^{a_{1}}(s)\in\Lambda_{a_i}^m$,
for all $i=2,\dots,k$ where
 \[\Lambda_a^m=\{s\in\Lambda^m,\,T_a(s)\le m-1\}\,.\]
\end{lemma}

\noindent {\it Proof}. The last property follows from the simple remark
that for $s\in\Lambda^m$ we have that  $\Theta^a(s)\neq s$ if and only if
$s\in\Lambda_a^m$. So, for the rest of the proof we suppose that all
transformations used verify the above property.

Let $\bar{s}^{(m)}$ be the sequence of
$\Lambda_m$ defined by $\bar{s}_1^{(m)}=1$ and
$\Delta\bar{s}_{k}^{(m)}= - \Delta\bar{s}_{k-1}^{(m)}$  for all $2\le k\le m$.
That is $\bar{s}^{(m)}\eqdef (0,1,0,1,\dots,0,1)$ if $m$ is odd and
$\bar{s}^{(m)}\eqdef (0,1,0,1,\dots,1,0)$ if $m$ is even.

 First we prove that
the statement of the lemma is equivalent to the following one: for
any sequence $s$ of $\Lambda_m$ such that $s\neq \bar{s}^{(m)}$, there
are integers $b_1,b_2,\dots,b_p\in\{0,1,2\}$ such that
\begin{equation}\label{sbar}
\bar{s}^{(m)}=\Theta^{b_p}\Theta^{b_{p-1}}\dots\Theta^{b_1}(s)\,.
\end{equation}
Indeed, suppose that the later property holds and let
$s'\in\Lambda_m$ such that $s'\neq s$. If $s'=\bar{s}^{(m)}$, then the
sequence $b_1,b_2,\dots,b_p$ satisfies the statement of the lemma.
If $s'\neq \bar{s}^{(m)}$, then let $c_1,\dots,c_l\in\{0,1,2\}$ such
that
\[\bar{s}^{(m)}=\Theta^{c_l}\Theta^{c_{l-1}}\dots\Theta^{c_1}(s')\,.\]
We notice that the transformations $\Theta^a$ are involutive, i.e. for all
$x\in\Lambda_m$,
\begin{equation}
\label{r} \Theta^a\Theta^a(x)=x.
\end{equation}
 Then we have
$\Theta^{c_1}\Theta^{c_{2}}\dots\Theta^{c_l}(\bar{s}^m)=s'$, so that
\[s'=\Theta^{c_1}\Theta^{c_{2}}\dots\Theta^{c_l}
\Theta^{b_{p}}\Theta^{b_{p-1}}\dots\Theta^{b_1}(s)\,,\] which
implies (\ref{ssp}). The fact that (\ref{ssp}) implies (\ref{sbar}) is
obvious.\\

Now we prove (\ref{sbar}) by induction in $m$.
It is not difficult to see that the result is
true for $m=1,2$ and $3$. Suppose that the result is true up to $m$ and
let $s\in\Lambda^{m+1}$ such that $s\neq\bar{s}^{(m+1)}$. For $j\le
m$, we  call $s^{(j)}$ the truncated sequence
$s^{(j)}=(s_0,s_1,\dots,s_j)\in\Lambda^j$.
From the hypothesis of induction, there exist
$b_1,b_2,\dots,b_p\in\{0,1,2\}$ such that
\begin{equation}\label{eqlem1}
\bar{s}^{(m)}=\Theta^{b_p}\Theta^{b_{p-1}}\dots\Theta^{b_1}(s^{(m)})\,
\end{equation}
where
\begin{equation}\label{eqlem2}
\mbox{$s^{(m)}\in\Lambda_{b_1}^m\;$  and $\;\Theta^{b_{i-1}}\Theta^{b_{i-2}}\dots
\Theta^{b_{1}}(s^{(m)})\in\Lambda_{b_i}^m$,$\;$ for all $i=2,\dots,p$.}
\end{equation}
Then, let us consider separately the case where $m$ is even and the case where $m$ is odd.

If $m$ is even and $\Delta s_m\Delta s_{m+1}=-1$, then we obtain directly that
$$\Theta^{b_p}\Theta^{b_{p-1}}\dots\Theta^{b_1}(s)=\bar{s}^{(m+1)}\,.$$
Indeed, from (\ref{eqlem2}), none of the transformations $\Theta^{b_{i-1}}\dots\Theta^{b_1}$,
$i=2,\dots,p$ affects the last step of $s$, so the identity follows from (\ref{eqlem1}).

If $m$ is even and $\Delta s_m\Delta s_{m+1}=1$, then from the hypothesis of induction
there exist $d_1,d_2,\dots,d_{r}\in\{0,1,2\}$ such that
\begin{equation}\label{eqlem3}\Theta^{d_r}\dots\Theta^{d_1}(s^{(m)})=
(\bar{s}^{(m-1)}, 2)
\end{equation}
which, from the above remark, may be chosen so that
\begin{equation}\label{eqlem4}
\mbox{$s^{(m)}\in\Lambda_{d_1}^m\;$  and $\;\Theta^{d_{i-1}}\Theta^{d_{i-2}}\dots
\Theta^{d_{1}}(s^{(m)})\in\Lambda_{d_i}^m$,$\;$ for all $i=2,\dots,r$.}
\end{equation}
Since from  (\ref{eqlem4}), none of the transformations $\Theta^{d_{i}}\dots\Theta^{d_1}$, $i=1,\dots,r$
affects the last step of $s$, it follows from (\ref{eqlem3}) that
\begin{equation}\label{eqlem5}\Theta^{d_r}\dots\Theta^{d_1}(s)=
(\bar{s}^{(m-1)}, 2, 3)\,.\end{equation}
Then by applying transformation $\Theta^2$, we obtain:
\begin{equation}\label{eqlem6}\Theta^2(\bar{s}^{(m-1)}, 2, 3)=(\bar{s}^{(m-1)}, 2, 1)\,.
\end{equation}
Hence, from (\ref{eqlem3}) and since none of the transformations
 $\Theta^{d_{r-i}}\dots\Theta^{d_r}$, $i=0,1,\dots,r-1$
affects the last step of $(\bar{s}^{(m-1)},2,1)$,  we have
\[\Theta^{d_1}\Theta^{d_2}\dots\Theta^{d_r}(\bar{s}^{(m-1)},2,1)=(s^{(m)},s_m-\Delta s_{m+1})\,.\]
Finally from (\ref{eqlem1}) and (\ref{eqlem2}), we have
\[\Theta^{b_p}\Theta^{b_{p-1}}\dots\Theta^{b_1}\Theta^{d_1}\dots\Theta^{d_r}\Theta^2\Theta^{d_r}\dots\Theta^{d_1}(s)=\overline{s}^{(m+1)}\]
and the induction hypothesis is true at the order $m+1$, when $m$ is even.\\

The proof when $m$ is odd is very similar and we will pass over some of the arguments
in this case. If $m$ is odd and $\Delta s_m\Delta s_{m+1}=-1$, then we obtain directly that
$$\Theta^{b_p}\Theta^{b_{p-1}}\dots\Theta^{b_1}(s)=\bar{s}^{(m+1)}\,.$$
If $m$ is odd and  $\Delta s_m\Delta s_{m+1}=1$ then from the hypothesis of induction,
there exist  $d_1,d_2,\dots,d_{r}\in\{0,1,2\}$ such that
\begin{equation}\label{eqlem8}
\Theta^{d_r}\dots\Theta^{d_1}(s^{(m)})=
(\bar{s}^{(m-1)}, -1)\end{equation}
and
\begin{equation}\label{eqlem9}
\mbox{$s^{(m)}\in\Lambda_{d_1}^m\;$  and $\;\Theta^{d_{i-1}}\Theta^{d_{i-2}}\dots
\Theta^{d_{1}}(s^{(m)})\in\Lambda_{d_i}^m$,$\;$ for all $i=2,\dots,r$.}
\end{equation}
Then it follows from (\ref{eqlem8}) and (\ref{eqlem9}) that
\begin{equation}\label{eqlem10}
\Theta^{d_r}\dots\Theta^{d_1}(s)=
(\bar{s}^{(m-1)}, -1, -2)\end{equation}
and by performing the transformation $\Theta^1\Theta^0\Theta^1=\Theta^{-1}$,
\begin{equation}\label{eqlem11}\Theta^0\Theta^1\Theta^0(\bar{s}^{(m-1)}, -1, -2)=(\bar{s}^{(m-1)}, -1, 0)\,.
\end{equation}
From (\ref{eqlem8}) and (\ref{eqlem9}), it follows that
$$\Theta^{d_1}\dots\Theta^{d_r}(\bar{s}^{(m-1)}, -1, 0)=(s^{(m)},s_m-\Delta s_{m+1}),$$
which finally gives from (\ref{eqlem1}) and (\ref{eqlem2}),
$$\Theta^{b_p}\Theta^{b_{p-1}}\dots\Theta^{b_1}\Theta^{d_1}\dots\Theta^{d_r}\Theta^0\Theta^1\Theta^0
\Theta^{d_r}\dots\Theta^{d_1}(s)=\bar{s}^{(m+1)}$$
and ends the proof of the lemma.~\QED

In the proof of Theorem \ref{prop2} for technical reasons we have to
consider two cases: $T_1<\infty$ a.s. and $\p(T_1 = \infty )>0$.
Lemma \ref{lem2} proves that in the first case $M$ is a divergent process.

\begin{lemma}\label{lem2} Any skip free process such that $T_1<\infty$ a.s.~and
$\Theta^a(M)\el M$ for $a=0$ and $1$ satisfies$:$
$$\lim_{n\rightarrow+\infty}\,[M]_n=+\infty\,,\;\;\;\mbox{a.s.}$$
\end{lemma}

\noindent {\it Proof}. Let us introduce the first exit time from
the interval $[-a,a]$:
\[\sigma_a(M)=\inf\{n:|M_n|=a\}\,,\]
where $a$ is any integer. Let us put
\[\Psi^a(M)=\left(\sum_{k=1}^n(\ind_{\{k\le\sigma_a\}}- \ind_{\{k>\sigma_a\}})\Delta M_k\right)_{n\geq 0}\,,\]
where $\Delta M_k=M_k-M_{k-1}$.
First we  observe that if $\Theta^a(M)\el M$ for $a=0$ and 1, then $\Psi^a(M)\el M$, for $a=0$ and 1.
This assertion is obvious for $a=0$ since $\sigma_0=T_0$. For $a=1$, it follows
from the almost sure identity:
\[\Psi^{a}(M)=\Theta^a(M)\ind_{\{T_a<T_{-a}\}}+\Theta^{-a}(M)\ind_{\{T_{-a}<T_a\}}\,.\]
and the equalities:
\begin{eqnarray*}
&&\{T_a(M)<T_{-a}(M)\}=\{T_a(\Theta^a(M))<T_{-a}(\Theta^a(M))\},\\\\
&&\{T_{-a}(M)<T_{a}(M)\}=\{T_{-a}(\Theta^{-a}(M))<T_{a}(\Theta^{-a}(M))\}\,.\end{eqnarray*}
Then from the almost sure inequality
\[\sigma_3(\Psi^1(M))\le \max\{T_1(M),T_{-1}(M)\}\,,\]
the fact that $T_1(M)<\infty$, $T_{-1}(M)<\infty$ a.s. and the
identity in law $\Psi^1(M)\el M$, we deduce that $\sigma_3(M)<+\infty$, a.s.
 It means, since $M$ is a symmetric
 process, that  $T_3(M)<\infty$ and $T_{-3}(M)<\infty$, a.s. Generalizing the above inequality, we obtain
 \[\sigma_{a+2}(\Psi^1(M))\le \max\{T_a(M),T_{-a}(M)\}\,.\] This
 gives in the same manner as before, that for each $a\geq 0$, $\sigma_a < \infty $ a.s.. From this it is
 not difficult to see that $\lim _{n\rightarrow \infty}[M]_n = +\infty$, $\p$-a.s..
\QED

The next lemma shows that in the case $\p(T_1=\infty)>0$ we can modify our
process $M$ by pasting to it an independent symmetric Bernoulli random
walk $S$ and reduce the case $\p(T_1=\infty)>0$ to the case $T_1 <
\infty$ a.s..

 We denote by $[M]_{\infty} = \lim _{k\rightarrow \infty}[M]_k$ which always
 exists since it is an increasing process and we put
$$T=\inf\{k\geq 0:[M]_k=[M]_\infty\}\,,$$
with $\inf \{\emptyset \} = +\infty$. We denote the extension of the
process $M$ by $X$ where for all $k\geq 0$
$$X_k = M_k \ind_{\{ k<T\}} + (M_T + S_{k-T})\ind_{\{k\geq T\}}.$$
Note that $X=M$, on the set $\{T=\infty\}$.
\begin{lemma}\label{lem3}
Let $M$ be a discrete skip free process which satisfies
$\Theta^a(M)\el M$ for some $a\in\mathbb{Z}$. Then $X$
also satisfies $\Theta^a(X)\el X$. Moreover, the $\sigma-$algebras generated by the respective
quadratic variations coincide, i.e.
$\sigma([M])=\sigma([X])$, $X$ is a divergent process $\p$-a.s. and $M=S^X_{[M]}$.
\end{lemma}\vspace{0.2cm}

\noindent {\it Proof}. We show that reflection property holds for
$X$. In this aim, we consider the
two processes $Y$ and $Z$ such that for all $k\geq 0$,
$$Y_k = \Theta^a(M)_k \ind_{\{ k<T\}} + (\Theta^a(M)_T - S_{k-T})\ind_{\{k\geq T\}},$$
$$Z_k = M_k \ind_{\{ k<T\}} + (M_T + \Theta^{a-M_T}(S)_{k-T})\ind_{\{k\geq T\}}.$$
We remark that
\begin{equation}\label{ast}
\Theta^a(X)=Y\ind_{\{T_a(Y)\leq T\}}+Z\ind_{\{T_a(Z)> T\}}
\end{equation}
and we write the same kind of decomposition for $X$:
\begin{equation}\label{ast2}
X=X\ind_{\{T_a(X)\leq T\}}+X\ind_{\{T_a(X)> T\}}\,.
\end{equation}
In view of  (\ref{ast}) and (\ref{ast2}), to obtain $X\el\Theta^a(X)$
it is sufficient to show that  for all bounded and measurable functional $F$,
$$\e [F(X)]=\e[F(Y)\ind_{\{T_a(Y)\leq T\}}]+\e[F(Z)\ind_{\{T_a(Z)> T\}}]\,.$$
Since reflection is a transformation which preserves the quadratic
variation of the process, the random time $T$ can be defined as a
functional of $Y$ as well as a functional of $Z$.
So we see that
the last equality is equivalent to $X\stackrel{\mathcal L}{=}Y$ and
$X\stackrel{\mathcal L}{=}Z$.
The first   equality in law follows from the fact that
$$ (M,S) \stackrel{\mathcal L}{=}(\Theta^a(M), -S)$$
which holds due to the reflection property of $M$ and $S$, and independency of
$M$ and $S$.  The second one holds since it can be reduced to the reflection property of $S$ itself,
by conditioning with respect to $M$.

Finally, the identity $M=S^X_{[M]}$  just follows from the construction of $X$.\QED

\noindent {\it Proof of Theorem $\ref{prop2}$}. Since both processes $M$ and $\Theta^a(M)$ have the same
quadratic variation, the identity in law of the statement is equivalent to:
for all $a=0,1,2$
\begin{equation*}
(M,[M])\el(\Theta^a(M),[M])\,.
\end{equation*}
Then we remark that the above equalities are equivalent to: for all $a=0,1,2$
\begin{equation*}
(S^M,[M])\el(S^{\Theta^a(M)},[M])\,.
\end{equation*}
Now it is crucial to observe the path by path equality: for each $a=0,1,2$
\[S^{\Theta^a(M)}=\Theta^a(S^M)\,,\]
from which we obtain
\begin{equation}\label{eq1p}
(S^M,[M])\el(\Theta^a(S^M),[M])\,.
\end{equation}
Hence,
\begin{equation}
\label{eq3p} \mathcal L( S^M | [M]) = \mathcal L(\Theta ^a(S^M) |
[M])\end{equation}
Fix $m\ge1$ and let $s,s'\in\Lambda^m$ with $s\neq s'$ be fixed. Consider the sequence of integers
$a_1,a_2,\dots,a_k\in\{0,1,2\}$ given in Lemma~\ref{lem1} such that
\begin{equation}\label{eq2p}
s=\Theta^{a_k}\Theta^{a_{k-1}}\dots\Theta^{a_1}(s')\,.
\end{equation}
Denote by $S^{M,m}$ the restricted path $(S_0,S_1,\dots,S_m)$.
Iterating (\ref{eq3p}), we may write for all $u\in \Lambda^m$ :
\[\p\left(S^{M,m}=u\,|\,[M]\right)=\p\left(\Theta^{a_1}\Theta^{a_2}\dots\Theta^{a_{k}}
(S^{M,m})=u\,|\,[M]\right)\,.\]
Applying (\ref{r}), we see that the right-hand side is equal to
\[\p\left(S^{M,m}=
\Theta^{a_k}\Theta^{a_{k-1}}\dots\Theta^{a_1}(u)\,|\,[M]\right)\,.\]
Take now $u=s'$ and use  (\ref{eq2p}), to obtain
\begin{equation}\label{00}
\p\left(S^{M,m}=s'\,|\,[M]\right)=\p\left(S^{M,m}=s\,|\,[M]\right)\,.
\end{equation}
If $T_1<\infty$ a.s. then from Lemma \ref{lem2} we can see that $M$ is
divergent and for all $m\geq 0$ $\p(S^{M,m}\in\Lambda^m)=1$.
 Then from (\ref{00}) the law of $S^{M,m}$ is uniform
over $\Lambda^m$ and  it coincides with the conditional law of $S^{M,m}$ given
$[M]$. Hence, $S^{M,m}$  is symmetric  Bernoulli
random walk on $[0, m]$ independent from  $[M]$.  Since this holds
 for all $m\ge0$, we conclude that $S^{M}$ is
a symmetric Bernoulli random walk which is independent of $[M]$. So from Definition \ref{def},
 $M$ is a divergent Ocone local martingale.

If $\p(T_1 = \infty)>0$, we consider the extension $X$ of the process
$M$ defined in Lemma \ref{lem3}. Then,  $X$  satisfies the hypotheses of Theorem \ref{prop2}. Moreover,
 from Lemma \ref{lem3}, $\p(T_1(X)<\infty)=1$. From what has just been proved
 $S^X$ is a symmetric Bernoulli random walk which is independent of
 $[X]$, and hence from $[M]$. This implies that $S^X$ and $[M]$ are
 independent. From  Lemma \ref{lem3} we have $M=S^X_{[M]}$, and, hence, the process $M$
is itself an Ocone martingale by Definition \ref{def}.\QED

\subsection{Counterexamples}\label{counterexample}

In this part, we give two examples of a discrete skip free process $M$ which satisfy
$M_0=0$, $\Theta^0(M)\el M$ and $\Theta^1(M)\el M$, but which are not a discrete
Ocone martingales.\vskip 0.5cm

\noindent {\bf Counterexample 1}: Let $(\epsilon_k)_{k\geq 1}$ be a sequence of
independent symmetric Bernoulli random variables. We put
$M_0=0,\, \Delta M_1= \epsilon_1,\,  \Delta M_2= \epsilon_2,\, \Delta
M_3= \epsilon_2$ and for $k>3$, $ \Delta M_k= \epsilon_k$. We
introduce also
$$M_n= \sum_{k=1}^n  \Delta M_k.$$
Since $[M]_n=n$ for all $n\geq 1$ and since $M$ is not Bernoulli
random walk, it can not be an Ocone martingale.

Let us verify that  $\Theta^a(M)\el M$ for $a\in \mathbb N \setminus
\{2\}$.
For $a=0$ we have reflection property since the $\epsilon _k$'s are
symmetric and independent.
For $a=1$ we consider four possible cases
related with the values of $(M_1, M_2, M_3)$. Let us put  $R_n
=\sum_{k=4}^n\epsilon_k$ for $n\geq 4.$

In fact, if $M_1=1, M_2=2, M_3=3$, we have
$ \Theta^1(M)=(0,1,0,-1,(-1-R_n)_{n\geq 4}))$\\
\noindent If $M_1=1, M_2=0, M_3=-1$, then
$ \Theta^1(M)=(0,1,2,3,(3-R_n)_{n\geq 4}))$\\
\noindent If $M_1=-1, M_2=0, M_3=1$, then
$ \Theta^1(M)=(0,-1,0,1,(1-R_n)_{n\geq 4}))$\\
\noindent If $M_1=-1, M_2=-2, M_3=-3$, then
$ \Theta^1(M)=(0,-1,-2,-3,\Theta^1(-3-(R_n)_{n\geq 4}))$

Similar presentation is valid for $M$:\\
if $M_1=1, M_2=2, M_3=3$, then
$ M =(0,1,2,3,(3+R_n)_{n\geq 4}))$,\\
if $M_1=1, M_2=0, M_3=-1$, then
$ M =(0,1,0,-1,(-1+R_n)_{n\geq 4}))$\\
if $M_1=-1, M_2=0, M_3=1$, then
$ M =(0,-1,0,1,(1+R_n)_{n\geq 4}))$\\
if $M_1=-1, M_2=-2, M_3=-3$, then
$ M =(0,-1,-2,-3,\Theta^1(-3+(R_n)_{n\geq 4}))$\\
To see that the laws of $\Theta^1(M)$ and $M$ are equal it is
convenient to pass to increments of corresponding processes.

If we take a pass with $M_1=1, M_2=2, M_3=3$, then $\Theta^2(M)$ of
such trajectory has a probability zero which is not the case for the
corresponding trajectory of $M$. So, $\Theta²(M)\nel M$.
For $a\geq 3$ we can write that
$$\Theta^3(M) = (M_1, M_2, M_3, \Theta^3((M_k)_{k\geq 4}))$$
and we conclude from symmetry of Bernoulli random walk.\vskip 0.5cm

\begin{center}
\includegraphics{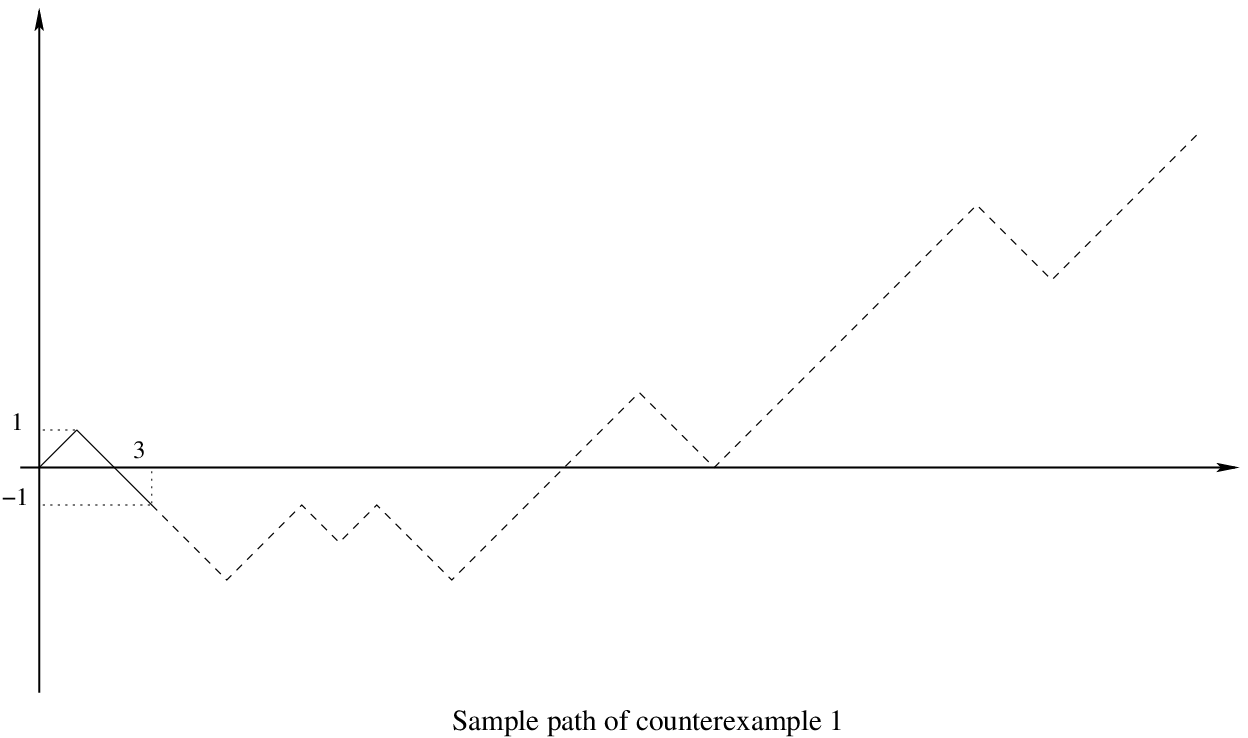}
\end{center}

\vspace*{.5in}

\noindent {\bf Counterexample 2}:
Let $(\varepsilon_k)_{k\ge0}$ be a sequence of independent $\{-1,+1\}$-valued symmetric Bernoulli
random variables. Set $k_n=\left\lfloor\frac{\ln (n+1)}{\ln 2}\right\rfloor-1$, where $\lfloor x\rfloor$
is the lower integer part of $x$ and let us consider the
following skip free process:
\[M_0=0\;\;\mbox{and for $n\ge1$,}\;\;M_n=\sum_{k=0}^{k_n} 2^k\varepsilon_k+(n-2^{k_n})\varepsilon_n\,.\]
Actually, $M$ is constructed as follows: $M_0=0$, $M_1=\varepsilon_0$ and for all $k\ge1$ and
 $n\in[2^{k},2^{k+1}-1]$, the increments $M_n-M_{n-1}$
have the sign of $\varepsilon_k$. In particular, the increments of $(M_n)$ are $-1$ or $1$ and since,
from the discussion at the beginning of section \ref{skipfree},  the only skip
free  local martingale with such increments is the Bernoulli random walk, it is clear that $M$ is not
an Ocone local martingale.\\

\begin{center}
\includegraphics{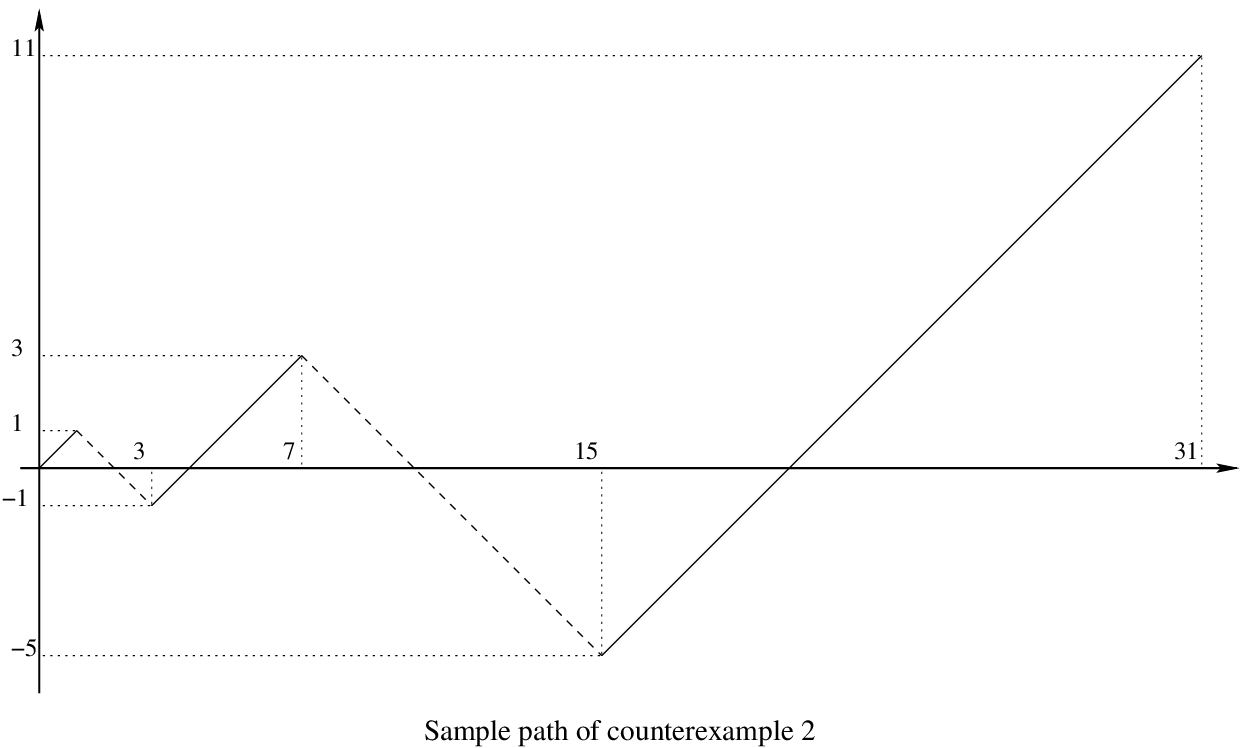}
\end{center}

\vspace*{.2in}

\noindent The  equality  $\Theta^0(M)\el M$  only means that $M$
is a symmetric process, which is straightforward  from its construction. Now let us check that $T_1<\infty$,
a.s.~and $\Theta^1(M)\el M$. Almost
surely on the set $\{M_1=-1\}$, there exists $k\ge0$ such that
$\varepsilon_i=-1$ for all $i\le k$ and  $\varepsilon_{k+1}=1$. The later assertion is equivalent to say that
for all integer $n\in(0,2^{k+1}-1]$, $M_n-M_{n-1}=-1$ and for all
$n\in[2^{k+1},2^{k+2}-1]$, $M_n-M_{n-1}=1$. It is then easy to check that
\[M_{2^{k+2}-1}=1\,.\]
So we have proved that $\{M_1=-1\}\subseteq\{T_{1}<\infty\}$, but since we also have
$\{M_1=1\}\subseteq\{T_{1}<\infty\}$, it follows that $\p(T_{1}<\infty)=1$.

Then we see from the construction of $(M_n)$
that almost surely, $T_1$ belongs to the set $\{2^j-1:j\ge1\}$ and that for $j\ge1$, conditionally to
$T_1=2^j-1$, $(M_n,\,n\le T_1)$ and $(M_{T_1+n},\,n\ge0)$ are independent. Moreover,
\[(M_{T_1+n},\,n\ge0)\el  (2-M_{T_1+n},\,n\ge0)\,,\]
so this proves that  $\Theta^1(M)\el M$.\\

Finally note that $0$ and $1$ are the only nonnegative levels at which the reflection principle
holds for the process $M$, i.e. $\Theta^a(M)\el M$ implies $a=0$ or $1$.  Indeed, at least it is clear
from the
construction of $M$ that the only times and levels at which the sign of its increments can change belong
to the set $\{2^j-1,\,j\ge0\}$,  i.e. if $a\ge0$ is such that $\Theta^a(M)\el M$, then necessarily
$a\in\{2^j-1,\,j\ge0\}$ and $T_a\in\{2^j-1,\,j\ge0\}$. But suppose that for $i\ge2$ we have $T_1=2^{i}-1$
and recall that all the increments  $M_{T_1+k+1}-M_{T_1+k}$ for all $k=0,1,\dots 2^i-1$ have the same sign. If
these increments are 1, then the process $M$  reaches the level $2^i-1$ at  time $T_1+2^i-2=2^{i+1}-3$ which
does not belong to the set $\{2^j-1,\,j\ge0\}$. So the sign of the increments of $M$ cannot change at this
time and the level $2^i-1$ cannot satisfy the identity in law $\Theta^{2^i-1}(M)\el M$.

\subsection{Continuous time lattice processes}\label{lattice}
As a preliminary result for the proof of Theorem \ref{main1}, we state an analogue of Theorem
\ref{prop2} for continuous time lattice processes.
We  say that $M=(M_t)_{t\geq 0}$ is a {\it continuous time lattice process} if $M_0=0$ and if it is a
pure jump c\`adl\`ag process whose jumps $\Delta M_t= M_t-M_{t-}$
verify : $|\Delta M_t|=\eta$, for some fixed real $\eta>0$. If we denote by
$(\tau_k)_{k\geq 1}$ the jump times of $M$, i.e. with $\tau_0=0$, for $k\ge1$,
\[\tau_k=\inf\{t>\tau_{k-1}:|M_t-M_{\tau_{k-1}}|=\eta\}\,,\]
with $\inf\{\emptyset\}=\tau_{k-1}$, then for all $t\ge0$ and $\p$-a.s.
$$M_t=\sum_{k=1}^{\infty}\Delta M_{\tau_k}\ind _{\{\tau_k\leq t\}}\,.$$
The quadratic variation  of $M$ is given by:
$$[M]_t=\sum_{k=1}^{\infty}(\Delta M_{\tau_k})^2\ind _{\{\tau_k\leq
  t\}}=\eta^2\sum_{k=1}^{\infty}\ind _{\{\tau_k\leq t\}}.$$
Note that $\tau_k$ admits the equivalent definition $\tau_k=\inf\{t\ge0:[M]_t=k\eta^2\}$.
We define the time changed discrete process $S^M$ by  $S^M=(M_{\tau_k})_{k\ge0}$ which has values
in the lattice $\eta\mathbb{Z}$. In particular, we have:
\begin{equation}\label{idlattice}
M_t=S^M_{\eta^{-2}[M]_t}\,,\quad t\ge0\,.
\end{equation}
We say that $M$ is a {\it continuous time lattice Ocone local martingale} if it can be written as
$M_t=S_{A_t}$, where $S$ is a symmetric Bernoulli random walk with values in the lattice $\eta\mathbb{Z}$
and $A$ is an increasing continuous time lattice process with values in $\mathbb{N}$ which is independent
of $S$. In the case where $M$ is divergent, $S$ coincide with $S^ M$ given in formula (\ref{idlattice}).
When $M$ is not divergent, $S$ is different from $S^ M$, namely if
$T= \inf \{k\geq 0: [S^ M]_k= [S^ M]_{\infty}\}$ then $S$ can be taken as:
$$S_k= S^M_k \ind_{\{ k\leq T\}} + (S^M_T + \tilde{S}_{T-k})\ind_{\{ k>T\}}\,,$$
where $\tilde{S}$ is a symmetric Bernoulli random walk which is independent from $S^M$.
In this case $S$ is independent from $[M]$. Therefore, when considering a continuous time lattice Ocone
local martingale $M$,  in  identity (\ref{idlattice}) we can and will suppose that $S^M$ is a symmetric Bernoulli
random walk with values in the lattice $\eta\mathbb{Z}$ and which is independent of $[M]$.
\vspace{0.2cm}

Recall the definitions (\ref{ta}) and  (\ref{transf}) of the hitting
time $T_{a}$ and transformations $\Theta^{a}$, respectively.

\begin{proposition}\label{c} Let $M$ be any continuous time lattice process
such that for all $k=0,1,2$, $$\Theta^{k\eta}(M)\el M\,,$$
then $M$ is a continuous time lattice Ocone local martingale.
If in addition  $T_{\eta}< \infty$ a.s.,
then  $S^M$ is a symmetric random walk on the lattice
$\eta\mathbb{Z}$ which is independent of $[M]$. Moreover, $M$ is a divergent local martingale
with respect to its own filtration.
\end{proposition}
\noindent {\it Proof}. Set $N=\eta^{-1}M$. We remark that for $k=1,2,3$,
$$\Theta^{k}(N)\el N.$$ Then following the proof of
Theorem \ref{prop2} along the lines for the continuous time process $N$, we obtain that $S^N$
conditionally to $[N]$ is Bernoulli random walk. Hence $S^{N}$ is a Bernoulli random walk which
is independent of $[N]$. Since $S^N=\eta^{-1}S^M$ and
$\eta^{-2}[N]=[M]$, we obtain that  $S^M$ is a symmetric Bernoulli
random walk on the lattice $\eta\mathbb{Z}$
which is independent of $[M]$. It means that it is local
martingale with respect to its own filtration. Finally, when $T_{\eta}<\infty$ a.s., $M$
 is a divergent local martingale since $N$ is so.\QED

\section{Proof of theorem 1}\label{proof}

Let $(\Omega ,\mathcal F , \mathbb F, \p)$ be the canonical space of continuous functions
with filtration $\mathbb F$ satisfying usual conditions.
Let $M$ be a continuous stochastic process which is defined on this space and
satisfying the assumptions of Theorem~\ref{main1}. Without loss of generality we
suppose that the sequence $(a_n)$ is decreasing.
\vskip 0.5cm

\noindent {\it Proof of Theorem $\ref{main1}$}. First of all we note
that since the map $x\rightarrow \Theta^x(\omega)$ is continuous on
$C(\R^+, \R)$, the hypothesis of this
theorem imply that $\Theta^0(M)\stackrel{\mathcal L}{=}M$, i.e. $M$ is
symmetric process.

Now, fix a positive integer $n$. We define the continuous lattice valued process
$M^n$ by using discretisation with respect to the space variable. In this aim, we
introduce the sequence
of stopping times $(\tau_k^n)_{k\geq 0}$  i.e. $\tau_0^n=0$ and for all $k\ge1$
$$\tau ^n_k= \inf\{t>\tau ^n_{k-1} \,: \, |M_t-M_{\tau ^n_{k-1}}| =a_n\}\,,$$
with $\inf\{\emptyset \}=\tau^n_{k-1}$. Then  $M^n=(M^n_t)_{t\geq0}$ is defined by:
$$M^n_t= \sum^{\infty}_{k=0} M_{\tau ^n_{k}}
\ind _{\{\tau^n_{k}\leq t < \tau ^n_{k+1}\}}\,.$$ We can easily check that  $M^n$
is a continuous time lattice process verifying the assumptions  of Proposition~\ref{c}.
Therefore according to this proposition, $[M^n]$ is a continuous time lattice Ocone
local martingale.

From the construction of $M^n$ we have  the almost sure inequality
\begin{equation}\label{bound}
\sup _{t\geq 0}|M_t -M^n_t|\leq a_n\,.
\end{equation}
Hence the sequence $(M^n)$ converges a.s.~uniformly on $[0,\infty)$ toward $M$.
The condition $$\sup_{n\geq 1}\sup_{t\geq 0}|\Delta M^n_t|\leq a_1$$
and (\ref{bound}) imply (cf. \cite{JS},Corollary IX.1.19, Corollary
VI.6.6) that $M$ is a local martingale and that
\begin{equation}\label{conv}
(M^n, [M^n])\stackrel{\mathcal L}{\rightarrow}(M,\langle M\rangle)\,.
\end{equation}
Since the properties (i) and (iii) given in introduction are equivalent, it is sufficient to verify that
for every deterministic function $h$ of the form $\sum_{j=1}^k\lambda_j\ind_{]t_{j-1}, t_{j}]}$
with $t_0=0 <t_1<\cdots t_k$ we have:
\begin{equation}\label{iii}
\e\left[\exp\left(i\int_0^\infty h(s)\,dM_s\right)\right]=
\e\left[\exp\left(-\frac12\int_0^\infty h^2(s)\,d\langle
M\rangle_s\right)\right]\,.
\end{equation}
From (\ref{conv}) we see that
$$\lim _{n\rightarrow \infty}\e\left[\exp\left(i\int_0^\infty h(s)\,dM_s^n\right)\right] =
\e\left[\exp\left(i\int_0^\infty h(s)\,dM_s\right)\right]\,.$$
Then in order to obtain (\ref{iii}), we will show by straightforward calculations that
\begin{equation}\label{calc}
\lim _{n\rightarrow \infty}\e\left[\exp\left(i\int_0^\infty h(s)\,dM_s^n\right)\right]=
\e\left[\exp\left(-\frac12\int_0^\infty h^2(s)\,d\langle
M\rangle_s\right)\right]\,.
\end{equation}
To prove (\ref{calc}) we first write
\[\e\left[\exp\left(i\int_0^\infty h(s)\,dM_s^n\right)\right]=
\int\e\left[\exp\left(i\int_0^\infty h(s)\,dM_s^n\right)\,|\,[M^n]=\omega\right]
\,dP_{[M^n]}(\omega)\,,\]
where $P_{[M^n]}$ is the law of $[M^n]$.
Then from Proposition 1 we have that
$$M^n \stackrel{\mathcal L}{=}a_nS_{a_n^{-2}[M^n]}$$
where $S$ is symmetric Bernoulli random walk independent from $[M^n]$.
Moreover,
$$\int_0^\infty h(s)\,dM_s^n = \sum _{j=1}^{k} \lambda _j \Delta M_{t_j}^n
\stackrel{\mathcal L}{=}a_n\sum _{j=1}^{k} \lambda _j \Delta S_{r_j}$$
where $\Delta M_{t_j}^n = M_{t_{j}}^n- M_{t_{j-1}}^n$, $\Delta S_{r_j}= S_{r_{j}} - S_{r_{j-1}}$ and
$r_j= a_n^{-2}[M^n]_{t_j}$, $1\le j\le k$.\\

Since $S$ and $[M^n]$ are independent and $\e\left[\exp (i a \Delta S_k)\right] = \cos (a)$ for all $a\in \R$, we have:
\begin{equation}\label{podivon}
\e\left[\exp\left(i\int_0^\infty h(s)\,dM_s^n\right)\,|\, [M^n]=\omega\right]=
\prod_{j=1}^k [\cos(\lambda _j a_n)]^{(u_j^n -u_{j-1}^n)}\,,
\end{equation}
where $u^n_j=\lfloor a_n^{-2}\omega_{t_j}\rfloor$, $j=0,1,\dots,k$ and $\lfloor x\rfloor$ is the lower integer part
of $x$. Moreover, it is not difficult to see that
\begin{equation}\label{hvatit}
\lim _{n\rightarrow \infty}\prod_{j=1}^k [\cos(\lambda _j a_n)]^{(u_j^n -u_{j-1}^n)} =
\exp\left(-\frac{1}{2}\sum_{j=1}^k \lambda_j^2(\omega_{t_j} - \omega_{t_{j-1}})\right)
\end{equation}
uniformly on compact sets of $\R_+^k$.
Then, the expression (\ref{podivon}) and the convergence relations (\ref{conv}), (\ref{hvatit}) imply (\ref{calc}). \QED

\section{Proof of Theorem \ref{main2}}\label{conj}
In what follows we assume, without loss of generality, that the process $M$ is divergent.
We begin with the following classical result of ergodic theory, a proof of which may be
found in \cite{dey}, Lemma 1.

Let $(\Omega ,\mathcal F , \mathbb F, \p)$ be canonical space of
continuous functions endowed by natural right-continuous filtration $\mathbb F =
(\mathcal F _t)_{t\geq 0}$ completed by negligible sets of $\mathcal
F= \bigvee_{t\geq 0} \mathcal F_t$.

\begin{lemma}\label{lemdey}
Let  $\Theta$ be a measurable transformation
of $\Omega$ to $\Omega$ which preserves $\p$. A random variable $X\in L^2(\Omega,{\cal F},\p)$ is
a.s.~invariant by $\Theta$ if and only if
\[\e (Z\cdot(Y\circ \Theta))=\e (Z\cdot Y)\,,\]
for all $Y\in L^2(\Omega,{\cal F},\p)$.
\end{lemma}\vskip 0.5cm

 Let $\Theta_n$, $n\ge1$ be a family of
transformations defined on canonical space of continuous  functions
 $(\Omega, \mathcal F, \mathbb F,\p)$.
Let ${\cal I}$ be the sub $\sigma$-algebra of the invariant events by all the transformations
$\Theta_n$, $n\ge1$, i.e.
\[{\cal I}=\{F\in{\cal F}:\ind_{F}\circ \Theta_n \stackrel{\mbox{a.s.}}{=}\ind_F,\;\mbox{for all $n\ge1$}\}.\]

 The following lemma extends Theorem 1 in \cite{dey}.
\begin{lemma}\label{lemcondind}
Let $M$ be a continuous  divergent local martingale defined on the filtered probability space
$(\Omega,\mathbb F, {\cal F},\p)$. Assume that the transformations $\Theta_n$ preserve the Wiener measure,
i.e.~if $B$ is the standard Brownian motion then for all $n\ge1$,
$B\circ \Theta_n\el B$.
The following assertions are equivalent:
\begin{itemize}
\item[$(j)$] For all $n\ge1$, $(B^M,\langle M\rangle)\el(\Theta_n(B^M),\langle M\rangle)$ have the same law.
\item[$(jj)$] $B^M$ and $\langle M\rangle$ are conditionally independent given the $\sigma$-field
$\mathcal I^M = (B^M)^{-1}({\cal I})$.
\end{itemize}
\end{lemma}
\noindent {\it Proof}. The proof almost follows from this of Theorem 1 in \cite{dey} along the lines.
 We first prove that $(j)$ implies $(jj)$.

Let $h,g$ two measurable functions $\Omega \rightarrow
\Omega$. Then (j) implies:
\begin{equation}\label{fu1}
\e\left(h(\langle M\rangle) g(B^M)\right)= \e\left(h(\langle M\rangle) g(B^M\circ
\Theta_n)\right) =\e\left(h(\langle M\rangle)( g(B^M)\circ \Theta_n)\right)
\end{equation}
We take conditional expectation with respect to $B^M$. For this we
denote by $f$ the following function:
$$\e \left(h(\langle M\rangle) | B^M\right)\stackrel{\mbox{a.s.}}{=}f(B^M)$$
Then (\ref{fu1}) implies that
\begin{equation}\label{fu2}
\e\left(f(B^M) g(B^M)\right)= \e\left(f(B^M)( g(B^M)\circ \Theta_n)\right)
\end{equation}
 Then according to Lemma \ref{lemdey} $f(B^M)$ is $\Theta_n$-invariant
 variable, i.e. it is measurable with respect to $\sigma$-algebra of
 $\Theta _n$-invariant sets $\mathcal I_n$. Since it holds for all $n\geq
 1$, $f(B^M)$ is measurable with respect to $\mathcal I=
 \cap_{n=1}^{\infty}\mathcal I_n$.
Moreover,\\\\
$\e\left(h(\langle M\rangle) g(B^M)| \mathcal I^M\right)=
\e\left(f(B^M) g(B^M) | \mathcal I\right)=$
$$\e \left(f(B^M) | \mathcal I\right)\e \left( g(B^M) | \mathcal I\right)=
\e (h(\langle M\rangle) | \mathcal I^M)\e(g(B^M) | \mathcal I^M)$$
and $(jj)$ is proved.

Now suppose that $(jj)$ is valid. Then
$$\e (h(\langle M\rangle) g(B^M))=
\e \left(\e \left(h(\langle M\rangle) | \mathcal
    I^M\right)\e \left(g(B^M) | \mathcal I^M \right)\right) $$

Moreover, since $B^M\circ \Theta_n$ and $\langle M \rangle $ are also conditionally independent for all
$n\ge1$, we have
 $$\e \left(h(\langle M\rangle) g(B^M\circ \Theta_n)| \mathcal I^M) \right)
=\e \left(h(\langle M\rangle) | \mathcal
 I^M\right)\e \left(g(B^M)\circ \Theta_n | \mathcal I^M\right) $$
Since every $ \mathcal I^M$-measurable random variable has the form $u(B^M)$,
where $u$ is ${\cal I}$-measurable,
$$\e\left(g(B^M)\circ \Theta_n)\,|\, \mathcal I^M\right) = \e\left(g(B^M)\,|\,
 \mathcal I^M\right)$$
and we obtain
$$\e \left(h(\langle M\rangle) g(B^M)\right)= \e\left(h(\langle M\rangle) g(B^M\circ
\Theta_n)\right)$$
which is $(j)$.\QED

\noindent\it Proof of Theorem \ref{main2}.\rm

If $(ii)$ holds then from Lemma \ref{lemcondind}, $B^M$ and $\langle M\rangle$ are
independent, so $(i)$ holds. Let us prove that $(i)$ implies $(ii)$. Suppose that $(ii)$ fails.
We show that $(i)$ fails, too.
Namely we show that one can construct a continuous martingale $M = B_{A}$, where $B$ is standard
Brownian motion and $A$ is non-decreasing
continuous adapted process, such that $M$ verify reflection properties of $(i)$
but it is not Ocone martingale.

 Let $X$ be a non trivial
$B^{-1}({\cal I}^{a})$-measurable bounded random variable. Call $({\cal F}_t^B)$
the natural filtration generated by $B$. Let $N_t=\e(X\,|\,{\cal F}_t^B) $ for all $t\geq 0$ and
$N=(N_t)_{t\geq 0}$. We remark that $N$ is a $({\cal
  F}^B_t)$-martingale invariant by all transformations
 $(\Theta^{a_n})$:
$$N \stackrel{\mathcal L}{=} N \circ \Theta^{a_n}.$$

Now, we can construct a finite non-constant stopping time
$T$ which is invariant by all the transformations $\Theta^{a_n}$ by
setting $\,T=\inf \{t\geq t_0 \,|\, N_t\in K\}$, where $t_0$ is large enough and
$K$ is a suitable Borel set. For instance we can choose $K$ such that
$\p(X\in K)\ge2/3$. Since $N_t\rightarrow X$ a.s. as $t\rightarrow
\infty$ we can find $t_0$ such that for $t\geq t_0$, $\p(N_t\in K)\ge1/2$.

Finally,  for $\alpha>0$, let us define
the following increasing process
\[A_t=\int_0^t\ind_{[0,T]}(s)+\alpha\ind_{]T,\infty[}(s)\,ds\,.\]
This process is not deterministic whenever $\alpha\neq1$ and since it is invariant by all the transformations
$\Theta^{a_n}$, one has $(B,A)\el(\Theta^{a_n}(B),A)$ for all $n\ge1$. The inverse of $A$ is given by
\[A_t^{-1}=\int_0^t\ind_{[0,T]}(s)+\alpha^{-1}\ind_{]T,\infty[}(s)\,ds\,,\]
so it is adapted and each $A_t$ is a $({\cal F}^B_t)$-stopping time.

Therefore $M=(M_t)_{t\geq 0}$ with $M_t=B_{A_t}$ is a continuous
divergent  $({\cal F}^B_{A_t})$-martingale satisfying $\Theta^{a_n}(M)\el M$, for all $n\ge1$. Moreover, $B^M=B$
and $\langle M\rangle=A$ are not independent by construction. Hence, $M$ can not be Ocone martingale with respect to the filtration
$({\cal F}^B_{A_t})_{t\geq 0}$ and it provides
a counterexample to the assertion $(i)$. So, we have proved that $(i)$ implies $(ii)$.\QED


\vspace*{.2in}

\begin{thebibliography}{99}

\bibitem{an} \sc D. Andr\'e: \rm  Solution directe du probl\`eme r\'esolu par M.~Bertrand.
{\it C.R. Acad. Sci. Paris}, {\bf 105}, 436-437, (1887).

\bibitem{oc} \sc D.L.~Ocone: \rm A symmetry characterization of conditionally
independent increment martingales. \it Barcelona
Seminar on Stochastic Analysis \rm 147--167, Progr. Probab., {\bf 32},
Birkh\"auser, Basel, 1993.

\bibitem{dey} \sc L.E.~Dubins, M.~\'Emery and M.~Yor: \rm
On the L\'evy transformation of Brownian motions and continuous martingales.
\it S\'eminaire de Probabilit\'es, \rm XXVII, 122--132, Lecture
Notes in Math., {\bf 1557}, Springer, Berlin, 1993.

\bibitem{ds} \sc L.E.~Dubins and M.~Smorodinsky: \rm The modified, discrete,
L\'evy-transformation is Bernoulli. {\it S\'eminaire de
Probabilit\'es}, XXVI, 157--161, Lecture Notes in Math., {\bf 1526},
Springer, Berlin, 1992.

\bibitem{JS} \sc J.~Jacod and  A.~Shiryaev: \it Limit Theorems for
Stochastic Processes, \rm Springer-Verlag Berlin, Heidelberg New York, 1987.

\bibitem{ma1} \sc M.~Malric: \rm Transformation de L\'evy et z\'eros du
mouvement brownien. {\it Probab. Theory Related Fields}, {\bf 101},
no. 2, 227--236, (1995).

\bibitem{ma2} \sc M.~Malric: \rm Densit\'e des z\'eros des transform\'es de L\'evy
it\'er\'es d'un mouvement brownien. {\it C. R. Math. Acad. Sci.
Paris}, {\bf 336}, no. 6, 499--504, (2003).

\bibitem{vy} \sc L.~Vostrikova and M.~Yor: \rm
Some invariance properties (of the laws) of Ocone's martingales. \it
S\'eminaire de Probabilit\'es, \rm XXXIV, 417--431, Lecture Notes in
Math., {\bf 1729}, Springer, Berlin, 2000.

\end{thebibliography}
\end{document}